\documentclass[11pt]{article} 
\usepackage{amsmath,amssymb,float,theorem} 

\makeatletter
\renewcommand{\@seccntformat}[1]{{\csname the#1\endcsname}.\hspace{5pt}}
\renewcommand{\section}{\@startsection{section}{2}{0em}{\baselineskip}
{\baselineskip}{\Large\bfseries}}
\renewcommand{\subsection}{\@startsection{subsection}{2}{0em}{\baselineskip}
{\baselineskip}{\large\bfseries}}%

\newtheorem{mythm}{Theorem}
\newenvironment{thm}{\begin{mythm}\hspace{-5pt}{\bf.}}{\end{mythm}}
\newtheorem{mylmm}{Lemma}
\newenvironment{lmm}{\begin{mylmm}\hspace{-5pt}{\bf.}}{\end{mylmm}}
\newtheorem{mycrl}{Corollary}
\newenvironment{crl}{\begin{mycrl}\hspace{-5pt}{\bf.}}{\end{mycrl}}
{\theorembodyfont{\normalfont}}

\def\BE#1{\begin{equation}\label{#1}}  
\def\EE{\end{equation}}

\def\ti#1{\widetilde{#1}} 
\def\e_ref#1{(\ref{#1})}

 \def\de{\delta}

\def\i{\infty}\def\p{\partial}\def\iv{^{-1}}
 \def\={\;=\;}  \def\+{\,+\,}  \def\a{\alpha}

\def\C{\mathbb C}  
  \def\hD{\widetilde D_w}
\def\F{\mathcal F}  \def\hF{\widetilde\F} 
\def\cH{\mathcal H}
\def\I{\mathcal I}   
\def\L{\mathfrak L} \def\M{\mathbf M} 
 \def\O{\textnormal O}
\def\cP{\mathcal P}\def\P{\cP}  
  
\def\Q{\mathbb Q} \def\S{\mathfrak S}
 \def\Z{\mathbb Z}

\begin{document}

\title{Some properties of hypergeometric series\\ 
        associated with mirror symmetry}
\author{Don Zagier and Aleksey Zinger}  
\date{\today}
\maketitle 

\begin{abstract}
We show that certain hypergeometric series used to formulate mirror symmetry
for Calabi-Yau hypersurfaces, in string theory and algebraic geometry, satisfy 
a number of interesting properties.  
Many of these properties are used in separate papers to verify the BCOV prediction 
for the genus one Gromov-Witten invariants of a quintic threefold and more 
generally to compute the genus one Gromov-Witten invariants of any 
Calabi-Yau projective hypersurface.  
\end{abstract}

\section{Introduction}
\label{mainthm_subs} 

\indent
An astounding prediction for the genus zero Gromov-Witten invariants of 
(counts of rational curves in) a quintic threefold was made in~\cite{CaDGP}.  
It was formulated in terms of the function~$\F$ defined in \e_ref{F0def} 
below and related objects.
This 1991 mirror symmetry prediction was mathematically verified about 
five years later.   
The 1993 mirror symmetry prediction of~\cite{BCOV} for the genus one 
Gromov-Witten invariants of a quintic threefold was recently verified in~\cite{Z1}.
A generalization of this prediction for a degree $n$ hypersurface $X_n$
in $\C P^{n-1}$, for an arbitrary~$n$, is proved in~\cite{Z2}; 
$X_5$ is a quintic threefold.
The proofs in these two papers make use of the properties of~$\F$ 
described by Theorems~\ref{HG_per}--\ref{HG_thm2} below.
Theorem~\ref{Phithm} explores related properties of~$\F$;
they appear to be of interest in their own right and may also be
of use in computation of higher genus Gromov-Witten invariants. 
Some further conjectural properties are stated in Section~\ref{conj_sec}.

We denote by
  $$\P\subset 1+x\Q(w)[[x]]$$
the subgroup of power series in $x$ with constant term~1 whose coefficients 
are rational functions in $w$ which are holomorphic at $w\!=\!0$. 
Thus, the evaluation map
  $$\cP\to 1+x\Q[[x]], \qquad F(w,x)\mapsto F(0,x)\,,$$
is well-defined. We define a map $\M:\P\to\P$ by
  $$\M F(w,x)\=\bigg\{1+\frac{x}{w}\frac{\p}{\p x}\bigg\}\frac{F(w,x)}{F(0,x)}\;.$$
Our first result says that the hypergeometric functions arising in 
the mirror symmetry predictions are periodic fixed points of the map~$\M$.

\begin{thm} \label{HG_per} Let $n$ be a positive integer and $\F\in\P$ 
the hypergeometric series
  \BE{F0def} \F(w,x)\=\sum_{d=0}^{\i}x^d
  \frac{\prod_{r=1}^{r=nd}(nw\!+\!r)}{\prod_{r=1}^{r=d}((w\!+\!r)^n\!-\!w^n)}\,.\EE
Then $\M^n\F=\F\,$.
\end{thm}
Note that we consider $n$ as fixed and therefore omit it from the notations.

If we now define further power series $\F_p\in\cP$ and $I_p\in 1+x\Q[[x]]$ 
for all $p\ge0$ by
 $$ \F_p(w,x)\=\M^p\F(w,x)\,,\qquad I_p(x)\=\F_p(0,x)\,,$$
so that $\F_{p+1}=(1+w^{-1}x\,d/dx)(\F_p/I_p)$, then Theorem~\ref{HG_per} 
says that $\F_{n+p}=\F_p$ and consequently $I_{n+p}=I_p$ for all $p\ge0$. 
The next result gives further properties of the functions $\{I_p\}_{p\in\Z/n\Z}\,$.
\begin{thm} \label{HG_thm1} The power series $I_p(x)$, $0\le p\le n-1$, satisfy
\begin{alignat}{1}
  \label{HG1e1} & I_0(x)\,I_1(x)\,\cdots\,I_{n-1}(x)\=(1-n^nx)\iv,\\
  \label{HG1e2} & I_0(x)^{n-1}I_1(x)^{n-2}\cdots I_{n-1}(x)^0\=(1-n^nx)^{-(n-1)/2}\,,\\
  \label{HG1e3} & I_p(x)\=I_{n-1-p} \qquad(0\le p\le n\!-\!1)\,.
\end{alignat}\end{thm}

We note that \e_ref{HG1e1} and the symmetry property~\e_ref{HG1e3} 
imply \e_ref{HG1e2}. 
However, \e_ref{HG1e2} is simpler to prove directly than~\e_ref{HG1e3} 
and will be verified together with~\e_ref{HG1e1} before we give 
the proof of~\e_ref{HG1e3}.

The power series $I_p$ describe the structure of $\F$ at $w\!=\!0$.
We will also describe some of its structure at $w\!=\!\i$.
We begin with the following observation, which will be proved in Subsection~\ref{HG2thm_subs}.

\begin{lmm} \label{reg_lmm} If $F\!\in\!\cP$ and $\M^kF=F$ for some $k>0$, 
then every coefficient of the power series $\,\log F(w,x)\in\Q(w)[[x]]$ is 
$\O(w)$ as $w\to\i\,$. \end{lmm}
Applying this lemma to $F=\F$, which satisfies its hypothesis 
by Theorem~\ref{HG_per}, we find that $\log \F(w,x)$ has an asymptotic 
expansion $\sum_{j=-1}^\i\mu_j(x)w^{-j}$ with $\mu_j(x)\in x\Q[[x]]$ 
for all $j\ge-1$ or equivalently, that $\F(w,x)$ itself has 
an asymptotic expansion 
\BE{asym} \F(w,x)\;\sim\;e^{\mu(x)w}\sum_{s=0}^{\i}\Phi_s(x)\,w^{-s} 
\qquad(w\to\i) \EE
for some power series 
$\mu=\mu_{-1},\;\Phi_0=e^{\mu_0},\;\Phi_1=\Phi_0\mu_1,\,\dots$ in $\Q[[x]]$.

\begin{thm} \label{HG_thm2}  
The first three coefficients $\mu(x)$, $\Phi_0(x)$, and $\Phi_1(x)$ 
in the expansion~$(\ref{asym})$ are given by
\BE{Th1.3} \mu(x)=\int_0^x\!\frac{L(u)-1}u\,du\,,\quad\; \Phi_0(x)=L(x)\,,
  \quad\; \Phi_1(x)=\frac{(n\!-\!2)(n\!+\!1)}{24n}\,\bigl(L(x)-L(x)^n\bigr)\,, \EE
where $L(x)$ denotes the power series $(1-n^nx)^{-1/n}\in\Z[[x]]\,$.
\end{thm}

The proof of this theorem in Subsection~\ref{HG2thm_subs} can be systematized 
and streamlined to obtain an algorithm for computing  every $\Phi_s$ 
by a differential recursion, which we now state. 
For integers $m\ge j\ge0$ (and for our fixed integer~$n$) we define $\cH_{m,j}=\cH_{m,j}(X)\in\Q[X]$ inductively by
\BE{Hdfn_e} \cH_{0,j}=\de_{0,j}, \qquad \cH_{m,j}=\cH_{m-1,j}+
(X-1)\bigg(X\frac{d}{dX}+\frac{m-j}{n}\bigg)\cH_{m-1,j-1}
\quad\text{for $m\ge1$} \EE
(with the convention that $\cH_{m-1,j-1}=0$ if $j=0$).
For example, for $0\le j\le 2$ we find
\BE{Hlow_e}\begin{split} \cH_{m,0}(X)&=1, \qquad
  \cH_{m,1}(X)\=\frac{1}{n}\binom{m}{2}(X-1), \\
  \cH_{m,2}(X)&\=\frac{1}{n^2}\binom{m}{3}\big((n+1)X-1\big)(X-1)
  \+\frac{3}{n^2}\binom{m}{4}(X-1)^2\,;\end{split}\EE
more generally, $\cH_{m,j}$ for fixed $j\ge1$ and varying $m$  
has the form $\sum_{k=1}^j\binom m{j+k}Q_{j,k}(X)$
with $Q_{j,k}\in\Z[n\iv,X]$ defined inductively~by 
$$Q_{0,k}=\de_{0,k}, \qquad 
Q_{j,k}=(X-1)(XQ_{j-1,k}'+(k\,Q_{j-1,k}+(k+j-1)\,Q_{j-1,k-1})/n 
\quad\text{for $j\ge1$.}$$
We now define differential operators $\L_k$ ($0\le k\le n$) on $\Q[[x]]$ by
\BE{Ldfn_e}
  \L_k\=\sum_{i=0}^k\bigg(\binom{n}{i}\cH_{n-i,k-i}(L^n) -(L^n-1)\sum_{r=1}^{k-i}
  \binom{n-r}{i}\frac{S_r(n)}{n^r}\,\cH_{n-i-r,k-i-r}(L^n)\bigg)D^i\,,\EE
where $D=x\,d/dx$ and $S_r(n)$ denotes 
the $r$th elementary symmetric function of $1,2,\ldots,n$
(a Stirling number of the first kind). 
Using~\e_ref{Hlow_e}, we find that the first two of these operators are
\begin{alignat}{1}
  \label{L1_e}  \L_1&\=nD-(L^n-1)=nLDL^{-1}\,,\\
  \label{L2_e}  \L_2&\=\binom{n}{2}D^2-\frac{3(n-1)}{2}(L^n-1)D
            \+\frac{n-1}{n}\bigg(\frac{(n-2)(n-11)}{24}L^n-1\bigg)(L^n-1)\,.
\end{alignat}

\begin{thm}\label{Phithm} $($i$)$ The power series $\Phi_s\in\Q[[x]]$, $s\ge0$, 
are determined by the first-order ODEs
\BE{PhiODE}
  \L_1(\Phi_s)\+\frac1L\,\L_2(\Phi_{s-1})\+\frac1{L^2}\L_3(\Phi_{s-2})
   \+\cdots\+\frac1{L^{n-1}}\,\L_n(\Phi_{s+1-n})\=0, \quad s\ge0,\EE
$($with the convention $\Phi_r=0$ for $r<0)$ together with the initial condition $\Phi_s(0)=\de_{0,s}\,$.  \\
$($ii$)$ For fixed $s$ and $n$, $\Phi_s(x)$ belongs to $L\,\Q[L]\,$. \\
$($iii$)$ For fixed $s$, $\Phi_s(x)$ belongs to $\Q(n)[L,L\iv,L^n]\,$.
\end{thm}

The meaning of part ($iii$) in Theorem~\ref{Phithm}
is that for each $s\ge0$ there exists 
$$\Psi_s\equiv\Psi_s(a,X,Y,Z)\in\Q(a)[X,Y,Z]$$
such that the function $\Phi_s(x)$ defined by~\e_ref{F0def} and~\e_ref{asym}
is given~by
$$\Phi_s(x)=\Psi_s\big(n,L(x),L(x)^{-1},L(x)^n\big).$$
In particular, ($iii$) neither implies nor is implied by ($ii$).

For example, from~\e_ref{PhiODE} for $s=0$ and $s=1$ together with equations \e_ref{L1_e} and~\e_ref{L2_e}
one finds the second and third identity in~\e_ref{Th1.3}, and continuing the same way one obtains
\begin{alignat*}{1}
\Phi_2 &= \frac{(n+1)^2(n-2)^2}{2\,(24n)^2}(L-2L^n+L^{2n-1})
=\Phi_1^2\big/2L\Phi_0\,,\\
\Phi_3 &= \frac{(n+1)(n-2)}{30\,(24n)^3}\,
\biggl\{(1003n^4-2366n^3+3759n^2-1676n-164)\,L^{3n-2}  \\
         &\qquad\qquad\qquad\qquad\quad\;
-\,72\,(n-1)(3n-1)(7n^2-9n+14)\,L^{2n-2} \\ 
         &\qquad\qquad\qquad\qquad\quad\;
+\,15\,(n+1)^2(n-2)^2\bigl(L^{2n-1}-L^n \bigr) \\
         &\qquad\qquad\qquad\qquad\quad\;
+\,72\,(n-1)(7n^3-17n^2+22n-24)\,L^{n-2} \\
         &\qquad\qquad\qquad\qquad\quad\;
+\,(5n^4+134n^3-447n^2+308n-556)\,L\biggr\}\,. 
\end{alignat*}
illustrating parts (ii) and (iii) of the theorem.  
These expressions, and the similar formulas obtained for $s\le7$, 
suggest that in fact $\Phi_s$ for $s$ fixed and $n$ varying is 
an element of $\Q[n,n\iv,L,L\iv,L^n]$, sharpening statement~(iii), 
but we do not know how to prove this.  Some further data and 
a further conjecture concerning the functions $\Phi_s(x)$ are given in Section~3.

\section{Proofs}  
\subsection{Preliminaries}

It will be convenient to introduce notations $D$ and $D_w$ for the first 
order differential operators $D=x\,\frac d{dx}$ and $D_w=D+w$ on $\Q(w)[[x]]$.  
(Here we think of $w$ as a parameter rather than a variable and 
write simply $\,\frac d{dx}\,$ instead of $\,\frac\p{\p x}$.)
The effect of $D_w$ on a power series $\sum c_d(w)x^d\in\Q(w)[[x]]$ 
is to multiply each $c_d(w)$ by $w+d$, so $D_w$ has an inverse operator 
$D_w\iv$ which replaces each $c_d(w)$ by $(w+d)\iv c_d(w)\,$. 
The operator $\,\M\,$ defined above can be written in terms of 
$D_w$ as $F(w,x)\mapsto w\iv D_w\bigl[F(w,x)/F(0,x)\bigr]\,$.

We remark that instead of working with the functions $\F_p(w,x)$, 
we could have worked with the functions $R_p(w,t)=e^{wt}\F_p(w,e^t)$, 
which are the objects that actually arise in the analysis of 
the mirror symmetry predictions for Gromov-Witten invariants. 
If we had done that, then the differential operator $D_w=w+x\,d/dx$ 
would have been replaced by the simpler differential operator $d/dt$, 
explaining why this operator plays such a ubiquitous role in our analysis.  
But it is easier, both in the calculations and for purposes of exposition, 
to work with power series over $\Q(w)$ in a single variable~$x$
rather than with objects in the less familiar space $e^{wt}\Q(w)[[e^t]]$.

The following lemma and its corollary are key to the proofs 
of the four theorems stated above.

\begin{lmm} 
\label{ode_lmm}
Suppose $c_0,\,\ldots,\,c_m,\,f,\,g,\,a$ are functions of $\,t\,$ 
(with $f$ not identically~$0$) satisfying 
  \BE{ode_lmm_e1}\begin{split}
  &c_m\,f^{(m)}\,+\,c_{m-1}\,f^{(m-1)}\,+\,\ldots\,+\,c_0\,f\=0\,,\\
  &c_m\,g^{(m)}\,+\,c_{m-1}\,g^{(m-1)}\,+\,\ldots\,+\,c_0\,g\=a\,, \end{split}\EE
where $f^{(k)}=d^k\!f/dt^k\,$.  Then the function $h:=(g/f)'$ satisfies
  \BE{ode_lmm_e2} \ti c_{m-1}\,h^{(m-1)}\,+\,\ti c_{m-2}\,h^{(m-2)}\,+\,\ldots\,+\,\ti c_0\,h\=a, \EE
where $\ti c_s(t)=\sum_{r=s+1}^m\binom{r}{s+1}\,c_r(t)\,f^{(r-1-s)}(t)\,$. 
\end{lmm}   
\noindent{\it Proof: } Using Leibnitz's rule and \e_ref{ode_lmm_e1}, we find
  $$ a\=\sum_{r=0}^mc_r\,\bigl(f\cdot g/f\bigr)^{(r)} \=\sum_{r=0}^mc_r\,\biggl(f^{(r)}g/f
  \,+\,\sum_{s=0}^{r-1}\,\binom{r}{s+1}f^{(r-1-s)}h^{(s)}\biggr)
 \=\sum_{s=0}^{m-1}\ti c_s\,h^{(s)}\,.$$

\begin{crl}\label{ode_cor}  
Suppose $F(w,x)\in\P$ satisfies
\BE{ode_cor_e1} 
\biggl(\sum_{r=0}^mC_r(x)\,D_w^{\,r}\biggr)F(w,x)\=A(w,x) \EE
for some power series $C_0(x),\,\dots,\,C_m(x)\in\Q[[x]]$ 
and $A(w,x)\in\Q(w)[[x]]$ with $A(0,x)\equiv0$.  
Then 
\BE{ode_cor_e2} \biggl(\sum_{s=0}^{m-1}\ti C_s(x)\,D_w^{\,s}\biggr)\M F(w,x)\=\frac1w\,A(w,x)\,, \EE
where $\,\ti C_s(x):=\sum_{r=s+1}^m\binom{r}{s+1}\,C_r(x)\,D^{r-1-s}F(0,x)\,$.   \end{crl}
\noindent{\it Proof: } Apply the lemma with $c_r(t)\!=\!C_r(e^t)$, $f(t)\!=\!F(0,e^t)$, $g(t)\!=\!e^{wt}F(w,e^t)$,
$a(t)\!=\!e^{wt}A(w,e^t)$, noting that then $h(t)=we^{wt}\M F(w,e^t)$.

\subsection{Proof of Theorem~\ref{HG_thm1}}  \label{Sec2.2}

For the proof of \e_ref{HG1e1} and \e_ref{HG1e2}, it is convenient to define $\F_p(w,x)$ also for $p=-1$. Set
  \BE{Ifuncdfn_e1a} \F_{-1}(w,x)\= w\,D_w\iv \F(w,x)\=\sum_{d=0}^{\i}x^d \frac{\prod_{r=0}^{r=nd-1}(nw\!+\!r)}
  {\prod_{r=1}^{r=d}\big((w\!+\!r)^n\!-\!w^n\big)}\;\in\;\P\,. \EE
We have $\F_{-1}(0,x)\!=\!1$ and $w\iv D_w\F_{-1}\!=\!\F$, 
so $\F_p=\M^{p+1}\F_{-1}$ for all $p\ge0$, justifying
the notation. It is straightforward to check that $\F_{-1}$ 
is a solution of the differential equation
\BE{HGODE_e1a} 
\biggl(D_w^n\,-\,x\,\prod_{j=0}^{n-1}\bigl(nD_w+j\bigr)\biggr)\,\F_{-1}
\=w^n\,\F_{-1}\,.\EE
This has the form of \e_ref{ode_cor_e1} with $F=\F_{-1}$, $A=w^n\F_{-1}$, $m=n$, and 
  \BE{init} C_n(x)\=1-n^nx\,, \quad C_r(x)\=-n^r\,S_{n-r}(n-1)\,x\qquad(0<r<n),\qquad C_0(x)\=0\,,\EE
where $S_{n-r}(n-1)$ as before denotes the $(n-r)$-th elementary symmetric function of 1,\,2,\,\dots,\,$n-1$.  
Applying Corollary~\ref{ode_cor} repeatedly, we obtain
\BE{ode-p} \sum_{s=0}^{n-1-p}C_s^{(p)}(x)\,D_w^s\F_p(w,x)\=w^{n-p-1}\,\F_{-1}(w,x)
\qquad (0\le p\le n-1), \EE
where $C_s^{(0)}(x)=C_{s+1}(x)$ with $C_r(x)$ as in \e_ref{init} 
and $C_s^{(p)}$ for $p>0$ is given inductively by
 \BE{inductive}   C_s^{(p)}\=\sum_{r=s+1}^{n-p}
 \binom{r}{s+1}\,C_r^{(p-1)}(x)\,D^{r-1-s}I_{p-1}(x)\,.\EE
In particular, by induction on $p$ we find that the first two coefficients 
in \e_ref{ode-p} are given by
\begin{alignat}{1}
\label{coeff1} 
C^{(p)}_{n-1-p} &\,\= (1-n^nx)\,\prod_{r=0}^{p-1}I_r(x)\,,\\
\label{coeff2} 
C^{(p)}_{n-2-p} &\,\= \biggl(-\frac{n^n(n-1)}2\,x
\,+\,(1-n^nx)\,\sum_{r=0}^{p-1}(n-r-1)\,
\frac{I_r'(x)}{I_r(x)}\biggr)\,\prod_{r=0}^{p-1}I_r(x)\,.
\end{alignat}
Equations \e_ref{ode-p} and \e_ref{coeff1} for  $p=n-1$ give
\BE{conseq1} 
(1-n^nx)\,\prod_{r=0}^{n-2}I_r(x)\,\F_{n-1}(w,x)\=\F_{-1}(w,x)\,.
\EE
Setting $w=0$ in this relation and using $\F_{-1}(0,x)=1$ 
gives equation \e_ref{HG1e1}. Then substituting \e_ref{HG1e1}
back into \e_ref{conseq1} gives $\F_{n-1}/I_{n-1}=\F_{-1}$ 
and hence, applying $w\iv D_w$ to both sides, $\F_n=\F$,
proving also part~(i) of Theorem~\ref{HG_thm1}. 
Similarly, taking $p=n-2$ in equations \e_ref{ode-p},  \e_ref{coeff1},  
and \e_ref{coeff2}  and then setting $w=0$ gives 
$$ \sum_{r=0}^{n-2} (n-r-1)\,\frac{I_r'(x)}{I_r(x)}\=\frac{n-1}2\,\frac{n^n\,x}{1-n^nx}\,,$$
and integrating this and exponentiating gives~\e_ref{HG1e2}.

Finally, we must prove the reflection symmetry~\e_ref{HG1e3}.  
For this purpose, it is useful to construct the power series $I_p$ in another way. 
Define a function $\hF_0\in\P$ by
\BE{newF0} \hF_0(w,x)\=\sum_{d=0}^{\i}x^d\,
\frac{\prod_{r=1}^{r=nd}(nw\!+\!r)}{\prod_{r=1}^{r=d}(w\!+\!r)^n}
\EE
and set $\hF_p(w,x)\=\M^p\hF_0(w,x)$ for all $p\ge0$.  
Since $\hF_0(w,x)$ is congruent to $\F(w,x)$ modulo~$w^n$,
we find by induction on $p$ that $\hF_p(w,x)$ is congruent to $\F_p(w,x)$ 
modulo~$w^{n-p}$ for all $0\le p\le n-1$ and
hence that $I_p(x)=\hF_p(0,x)$ in this range. 
We now argue as above, using $\hF_0$ instead of $\F_{-1}$. 
This function satisfies the differential equation
$$\biggl(D_w^{n-1}\,-\,nx\prod_{j=1}^{n-1}\bigl(nD_w+j\bigr)\biggr)\,\hF_0
\=w^{n-1}\,.$$
Applying Corollary~\ref{ode_cor} repeatedly, we obtain
$$ \sum_{s=0}^{n-1-p}\ti C_s^{(p)}(x)\,D_w^s\hF_p(w,x)\=w^{n-p-1}  $$
for $0\le p\le n-1$, where the coefficients $\ti C_s^{(p)}(x)\in\Q[[x]]$ 
can be calculated recursively, the top one being given by 
$$\ti C_{n-1-p}^{(p)}(x)=(1-n^nx)I_0(x)\cdots I_{p-1}(x).$$ 
Specializing to $p=n-1$ and using~\e_ref{HG1e1}, 
we find that $\hF_{n-1}(w,x)=I_{n-1}(x)$ is independent of~$w$.  
Now by downwards induction on~$p$, 
using the equation $\hF_p=I_pwD_w\iv\hF_{p+1}$, 
we can ``reconstruct'' all of the power series $\hF_p(w,x)$
($n-1\ge p\ge0$) from their special values $I_p(x)=\hF_p(0,x)$ at $w=0$, 
obtaining in particular the formula
$$ w^{1-n}\,\hF_0(w,x)
\=I_0\,D_w\iv\,I_1\,D_w\iv\,\cdots I_{n-2}\,D_w\iv\, I_{n-1} $$
for the initial series $\hF_0$.  
Comparing the coefficients of $x^d$ on both sides of this equation, we find
$$ \frac{n\iv\prod_{r=0}^{nd}(nw+r)}{[w(w+1)\cdots(w+d)]^n}\=  \sum_{\genfrac{}{}{0pt}{}{d_0,\dots,d_{n-1}\ge0}{d_0+\cdots+d_{n-1}=d}}
\frac{c_0(d_0)\,\cdots\,c_{n-1}(d_{n-1})}
{(w+d_1+\cdots+d_{n-1})(w+d_2+\cdots+d_{n-1})\cdots(w+d_{n-1})}$$
for all $d\ge0$, where $c_p(d)$ denotes the coefficient of $x^d$ 
in $I_p(x)$. 
Splitting up the sum on the right into the subsum over $n$-tuples
$(d_0,\dots,d_{n-1})$ with $\,\max\{d_r\}\le d-1$ and 
the sum over the $n$-tuples which are permutations of $(d,0,\dots,0)$, 
and using that $c_p(0)=1$ for all $p$, we can rewrite this equation as
$$\sum_{p=0}^{n-1}\frac{c_p(d)}{w^{n-p-1}(w+d)^p}
\=\frac{\prod_{r=0}^{nd}(nw+r)}{n\prod_{r=0}^d(w+r)^n} \,-\, 
\sum_{\genfrac{}{}{0pt}{}{0\le d_0,\dots,d_{n-1}<d}{d_0+\cdots+d_{n-1}=d}}
\frac{c_0(d_0)\,\cdots\,c_{n-1}(d_{n-1})}{(w+d_1+\cdots+d_{n-1})\cdots(w+d_{n-1})}$$
Now suppose by induction that $c_p(d')=c_{n-p-1}(d')$ for all $d'<d$ and
all $0\le p\le n-1$. 
(Notice that this is true for $d'=0$ because $c_p(0)=I_p(0)=1$ for all $p$, 
providing the starting point for the induction.)  
Then both terms on the right are $(-1)^{n-1}$-invariant under the map $w\to-w-d$, 
as one sees for the second term by making the renumbering $d_r\to d_{n-1-r}$.  
It follows that the left-hand side has the same invariance and hence that $c_p(d)=c_{n-1-p}(d)$ for all $0\le p\le n-1$, 
completing the inductive proof of the desired symmetry $I_{n-1-p}=I_p$.

\subsection{Proof of Theorem~\ref{HG_thm2}}
\label{HG2thm_subs}

We now turn to the expansion of $\F(w,x)$ near $w=\i$.  We first prove Lemma~\ref{reg_lmm}, which said
that any periodic fixed point of the map $\M:\P\to\P$ has a logarithm which belongs to $w\,\Q[[x,w\iv]]$.

\noindent{\it Proof of Lemma~\ref{reg_lmm}:}
The effect of $\M$ on logarithms is given~by $\M\bigl(e^{H(w,x)}\bigr)=e^{H^*(w,x)}\,$, where
  \BE{Mlog_e}  H^*(w,x)\=H(w,x)-H(0,x)+
              \log\biggl(1+\frac{DH(w,x)-DH(0,x)}w\biggr); \EE
here, as before, $D$ denotes $x\frac{\p}{\p x}$.  
Suppose that $H(w,x):=\log F(w,x)$ is not $\O(w)$, and let $e$ be the
smallest integer such that the coefficient of $x^e$ in $H(w,x)$ 
is not $\O(w)$ as $w\to\i$. Then
\BE{Hexp_e} H(w,x)\=Cx^ew^N+x \O_w(w)+x^e \O_w(w^{N-1})+ \O(x^{e+1})  
\qquad(w\to\i) \EE
for some $C\ne0$ and $N\ge2$, where $\O_w(w^\nu)\,$ denotes a polynomial in~$x$ 
with coefficients that grow at most like $w^\nu$ as $w\to\i$ and $\O(x^{e+1})\,$
denotes an element of $x^{e+1}\,\Q(w)[[x]]$. 
From~\e_ref{Mlog_e} and~\e_ref{Hexp_e},
 $$H^*(w,x)\=H(w,x)+Cex^ew^{N-1}+x\O_w(1)+x^e\O_w(w^{N-2})+\O(x^{e+1}).$$
This has the same form as \e_ref{Hexp_e} with the same $C$, $e$, and $N$.  
Iterating, we find that
$$\log\big(\M^kF(w,x)\big)\=H(w,x)+kCex^ew^{N-1}+
      x\O_w(1)+x^e\O_w(w^{N-2})+\O(x^{e+1}),$$
and this contradicts the assumption that $\M^kF=F$, since $C\ne0$
and $N\ge2$.

As already mentioned in the introduction, Lemma~\ref{reg_lmm} 
together with Theorem~\ref{HG_per} implies 
that $\F(w,x)$ has an asymptotic expansion of the form \e_ref{asym}.  
From the proof of the lemma, we see that each $\F_p(w,x)=\M^p\F(w,x)$
has an asymptotic expansion
  \BE{asym_p} \F_p(w,x)\;\sim\;e^{\mu(x)w}\sum_{s=0}^{\i}\Phi_{p,s}(x)\,w^{-s} \qquad(w\to\i) \EE
of the same form, with the same function $\mu(x)$ in the exponent.  
The equation $\F_{p+1}\!=\!\M \F_p$ gives
  \BE{expind_e}  \Phi_{0,s}=\Phi_s,\quad \Phi_{p+1,s}=\frac{1+\mu'}{I_p}\,\Phi_{p,s}\,+\, \begin{cases}
    \bigl(\Phi_{p,s-1}/I_p\bigr)'&\hbox{if}~s\ge1,\\ 
 \qquad\quad 0&\hbox{otherwise}, \end{cases}\EE
where $f'$ denotes $Df=x\,df/dx\,$.  
We want to solve these equations by induction on $p$ for small~$s$.

Before doing this, we begin with the following observation.  
Let $L(x)=(1-n^nx)^{-1/n}$ as in Theorem~\ref{HG_thm2}.
Then  \e_ref{HG1e1} says that the product of the functions $I_p(x)/L(x)$ ($p\in\Z/n\Z$) equals 1, 
so if we define
\BE{Hp} H_p(x)\=\frac{L(x)^p}{I_0(x)\cdots I_{p-1}(x)} \qquad(p\ge0), \EE
then we have the properties
\BE{HpProperties} H_0=1,\quad H_p/H_{p+1}=I_p/L, \quad H_1H_2\cdots H_n=1,
\quad H_{p+n}=H_p, \quad H_{n-p}=H_p\iv, \EE
where the last equality is originally true for $0\le p\le n$ but then, 
in view of the periodicity of $\{H_p\}$, holds for any $p\in\Z/n\Z$.
A number of identities below are simpler to state in terms of 
the functions $H_p(x)$ than in terms of the original functions~$I_p(x)$.

The case $s=0$ of \e_ref{expind_e} gives by induction  the formula
$\Phi_{p,0}=(1+\mu')^p/I_0\cdots I_{p-1}$.  
Combining this with the formulas $\F_n=\F$ and~\e_ref{HG1e1},
we obtain $(1+\mu')^n=L^n\,$, from which  
the first equation in~\e_ref{Th1.3} follows since $\mu(x)$ is a
power series in~$x$ with no constant term. This also gives us the formula 
$$ \Phi_{p,0}(x) \= H_p(x)\,\Phi_0(x)\qquad\text{for all $p\ge0$},  $$  
with $H_p$ as in \e_ref{Hp}.  
Now substituting this into the case $s=1$ of \e_ref{expind_e} we find inductively
$$\Phi_{p,1}(x)\=H_p(x)\,\biggl(\,\Phi_1(x)\,+\,p\,\frac{\Phi_0'-L'}L\,+\,
\frac{\Phi_0}{L}\sum_{r=1}^p\frac{H_r'}{H_r}\biggr)\qquad\text{for all $p\ge0$}.$$  
Setting $p=n$ in this relation and using the third and fourth of equations \e_ref{HpProperties} and $\F_n=\F$, we deduce that 
$\Phi_0=L$, which is the second assertion of Theorem~\ref{HG_thm2}.
At the same time we can refine the last two equations to
\BE{Phi0}  \Phi_{p,0} \= H_p\,L,\qquad
    \Phi_{p,1}\=H_p\,\biggl(\Phi_1 \,+\,\sum_{r=1}^p\frac{H_r'}{H_r}\biggr)\qquad\quad(p\ge0).\EE 

The proof of the third identity in \e_ref{Th1.3} 
is similar, but the calculations are more complicated. 
The case $s=2$ of \e_ref{expind_e} gives by induction the formula
$$\Phi_{p,2}\= H_p\,\biggl(\Phi_2\,+\,p\,\Big(\frac{\Phi_1}L\Big)' 
\,+\,\biggl(\sum_{r=1}^p\frac{H_r'}{H_r}\biggr)\frac{\Phi_1}L 
\,+\,\frac{1}{L}
\sum_{s=2}^p\sum_{r=1}^{s-1}\frac{H_r'}{H_r}\frac{H_s'}{H_s}
\,+\,\biggl(\frac{1}{L}
\sum_{r=1}^{p-1}(p\!-\!r)\frac{H_r'}{H_r}\biggr)'\biggr)$$
for all $p\ge0$. 
Taking $p=n$, observing that 
$$\sum_{s=2}^n\sum_{r=1}^{s-1}\frac{H_r'}{H_r}\frac{H_s'}{H_s}
\equiv\frac{1}{2}\Bigg(\bigg(\sum_{p=1}^n\frac{H_p'}{H_p}\bigg)^2
-\sum_{p=1}^n\bigg(\frac{H_p'}{H_p}\bigg)^2\Bigg)
=-\frac{1}{2}\sum_{p=1}^n\bigg(\frac{H_p'}{H_p}\bigg)^2$$
by the third equation in~\e_ref{HpProperties}, 
and using $\F_n=\F$, we find that 
$$n\Big(\frac{\Phi_1}L\Big)'
\=\frac{1}{2L}\sum_{p=1}^n\bigg(\frac{H_p'}{H_p}\bigg)^2
\,+\,\biggl(\frac{1}{L}\sum_{p=0}^{n-1}p\,\frac{H_p'}{H_p}\biggr)'
 \=-\frac{(n+1)(n-2)}{24}\,\bigl(L^{n-1}\bigr)'\,,$$
the last equation being Lemma~\ref{3rdDesc} below.  
Integrating and using $\Phi_1(0)=0$ gives the last identity in~\e_ref{Th1.3}.

\begin{lmm}\label{3rdDesc} 
The functions $\{H_p(x)\}_{p\in\Z/n\Z}$ satisfy
\begin{equation}\label{3rdDesc_e}
\frac1{2L}\sum_{p\;(\text{\rm mod $n$})}\biggl(\frac{H_p'}{H_p}\biggr)^2 
\= -\biggl(\frac{(n+1)(n-2)}{24}\,L^{n-1}\,+\,
  \frac1L\sum_{p=0}^{n-1}\,p\,\frac{H_p'}{H_p}\biggr)' \;. 
\end{equation}
\end{lmm}

The proof consists of expressing the left-hand side of \e_ref{3rdDesc_e}
in terms of the functions $I_0,I_1,\ldots,I_{n-1}$ and their derivatives,
getting rid of all square terms via the product rule, and 
then eliminating $I_{n-1}$, $I_{n-2}$, and~$I_{n-3}$.
The last elimination is achieved by computing the coefficients $C_p^{(n-3-p)}$ inductively by~\e_ref{ode-p}, starting with $$C_{n-3}^{(0)}\=-n^{n-2}\,S_2(n-1)\,x=-\frac{(n-1)(n-2)(3n-1)}{24}L'/L^{n+1}\,,$$
and then setting $p=n-3$, exactly as we did with $C_p^{(n-1-p)}$ 
and $C_p^{(n-2-p)}$ in Subsection~\ref{Sec2.2} to prove eqs.~\e_ref{HG1e1} and~\e_ref{HG1e2}.   
At this stage, all terms involving products of two functions $I_p$ cancel, 
and the resulting expression can be integrated. 
We omit the details, which are somewhat tedious, since the last identity in~\e_ref{Th1.3} also follows easily from Theorem~\ref{Phithm}.

\subsection{Proof of Theorem~\ref{Phithm}}
\label{Phithm_subs}
We set $X=L^n$ and $Y=(L^n-1)/n$. Note that 
  \BE{Phipf_e0} D(\mu)=L-1, \qquad D(L)=LY, \qquad D(X)=X^2-X, \qquad D(Y)=XY.\EE
The first identity implies that $D_w\,e^{\mu w}=e^{\mu w}\,\hD$, where $\hD=D+Lw$. 
By induction on~$k$, the powers of the differential operator $\hD$ are given by
  \BE{Dind_e1}\begin{split}  
  \hD^k &\=\sum_{m=0}^k\binom{k}{m}\,\hD^m(1)\,D^{k-m}\\
  &=D^k+k\,Lw\,D^{k-1}+\frac{k(k-1)}2\,\big((Lw)^2+Y(Lw)\big)\,D^{k-2}+\ldots\,.  
  \end{split}\EE
A second induction gives the formula
  \BE{Dind_e2}\hD^m(1) \=\sum_{j=0}^m \cH_{m,j}(X)\,(Lw)^{m-j},\EE
with $\cH_{m,j}\in\Z[X,Y]\subseteq\Q[X]$ given by~\e_ref{Hdfn_e}.

The function $\F(w,x)$ satisfies the ODE 
$$\bigg(D_w^n\,-\,w^n\,-\,x\,\prod_{j=1}^n\bigl(nD_w+j\bigr)\bigg)\F\=0\,.$$
Since $D_w\,e^{\mu w}=e^{\mu w}\,\hD$, the function $\hF(w,x)=e^{-\mu(x)w}\F(w,x)$ 
satisfies
the differential equation $\L\hF=0$, where $\L$ is the differential operator
 \begin{equation*}\begin{split}
  \L &\= L^n\,\biggl(\hD^n\,-\,w^n\,-\,x\,\prod_{j=1}^n\bigl(n\hD+j\bigr)\biggr) \\ 
  &\= \hD^n\,-\,(Lw)^n\,-\,(L^n-1)\,\sum_{r=1}^n\frac{S_r(n)}{n^r}\,\hD^{n-r}\,. 
 \end{split}\end{equation*}
Using \e_ref{Dind_e1} and~\e_ref{Dind_e2}, we can expand $\L$ as
$\L=\sum_{k=1}^n(Lw)^{n-k}\L_k$, with  $\L_k$ defined by~\e_ref{Ldfn_e}.
Combining the differential equation $\L\hF=0$ with the asymptotic expansion 
$\hF(w,x)\sim\sum_{s\ge0}\Phi_s(x)w^{-s}$ for large~$w$, we obtain~\e_ref{PhiODE}.

We will next use \e_ref{PhiODE} to prove by induction that $\Phi_s$ belongs 
to $L\Q[L]$.  
Since $\L_1(L\Q[L])=L^2Y\Q[L]$, it suffices to show that
 \BE{Lstat_e1}\L_k(L\Q[L])\subseteq L^{k+1}Y\,\Q[L] \qquad(2\le k\le n).\EE
Let $\I\subset\Q[L]$ be the ideal generated by $XY$.
Since $D$ and $Y$ commute modulo $\I$ by~\e_ref{Phipf_e0} and 
since $(D-rY)L^r=0$, we have
 \begin{equation*}
  (D-Y)(D-2Y)\ldots(D-kY)L^r\;\in \;
  \begin{cases} L^r\,\I &\hbox{if}~1\le r\le k\,,\\ L^rY\,\Q[L] &\hbox{if}~r\ge k+1\,.
 \end{cases}\end{equation*}
Therefore \e_ref{Lstat_e1} is a consequence of the following lemma.
\begin{lmm}\label{Lmod_lmm}
For all $k>1$,
  $$\L_k\;\equiv\;\binom{n}{k}(D-Y)(D-2Y)\cdots(D-kY)\;\pmod\I\,.$$
\end{lmm}
\noindent{\it Proof:} The recursion~\e_ref{Hdfn_e} for $H_{m,j}$ shows that 
$H_{m,j}\equiv h_{m,j}\,Y^j\pmod\I$, where $h_{m,j}\in\Z$ is given recursively by
  \BE{hdfn_e} h_{0,j}=\de_{0,j}, \qquad 
   h_{m,j}=h_{m-1,j}+(m-j)h_{m-1,j-1} ~~~\forall~m\ge1\,\EE
(with $h_{m-1,j-1}=0$ for $j=0$). 
Thus $h_{m,j}=\S_m^{(m-j)}$, where $\S_m^{(k)}$ denotes a Stirling number 
of the second kind (the number of ways of partitioning a set of $m$ 
elements into $k$ non-empty subsets).
We also note $(1-L^n)n^{-r}\equiv(-1)^rY^r\pmod\I$ for all $r\ge1$.  
Combining these facts with~\e_ref{Ldfn_e}, we find that
  $$\L_k\;\equiv\;
  \sum_{i=0}^k\biggr(\sum_{r=0}^{k-i}(-1)^r\binom{n-r}i\,S_r(n)\,
  \S_{n-r-i}^{(n-k)}\biggr)\,Y^{k-i}D^i\;\pmod\I\,.$$
The desired congruence for $\L_k$ now follows from the generating series calculation
 \begin{alignat*}{1}
 &\sum_{i=0}^k\biggr(\sum_{r=0}^{k-i}(-1)^r
  \binom{n-r}{i}\,S_r(n)\,\S_{n-r-i}^{(n-k)}\biggr)\,t^i \\
 &\= \sum_{i=0}^n\biggr(\sum_{r=0}^{n-i}(-1)^r\binom{n-r}i\,S_r(n)\,
  \biggl[\frac1{(n-k)!}\sum_{j=0}^{n-k}(-1)^{n-k-j} 
  \binom{n-k}jj^{n-r-i}\biggr]\biggr)\,t^i \\
 &\=\frac1{(n-k)!}\sum_{j=0}^{n-k}(-1)^{n-k-j}
  \binom{n-k}{j}\sum_{r=0}^n(-1)^rS_r(n)
  \sum_{i=0}^{n-r}\binom{n-r}ij^{n-r-i}t^i \\
 &\=\frac1{(n-k)!}\sum_{j=0}^{n-k}(-1)^{n-k-j}\binom{n-k}j
  \sum_{r=0}^n(-1)^rS_r(n)(j+t)^{n-r} \\
 &\=\frac{n!}{(n-k)!}\sum_{j=0}^{n-k}(-1)^{n-k-j}\binom{n-k}j\binom{j+t-1}n \\
 &\=\frac{n!}{(n-k)!}\binom{t-1}k \=\;\binom nk\,(t-1)(t-2)\cdots(t-k)\,,
 \end{alignat*} 
where the first equality follows from the well-known fact that 
the expression in square brackets equals $\S_{n-r-i}^{(n-k)}$ if $i+r\le k$ 
and 0 for $i+r>k$ and the second-to-last equality is obtained
by expanding $(1+u)^{t-1}((1+u)-1)^{t-1}$ by the binomial theorem 
and equating coefficients of $t^n$.

This completes the proof of part (ii) of Theorem~\ref{Phithm}.
Part~(iii) of Theorem~\ref{Phithm} follows from 
the differential equation~\e_ref{PhiODE} by induction on~$s$.

\section{Further discussion of the large $w$ expansion of $\F(w,x)$}
\label{conj_sec}
In this final section we give some further information and conjectures about 
the power series $\Phi_s(x)$ defined by equation~\e_ref{asym}.  
We begin by giving the numerical values for $n\le5$ and $s\le4$.  
For this purpose it is convenient to divide $\Phi_s/L$ by $((n-2)(n+1)/24n)^s/s!$ 
and write the result as the sum of $(1-L^{n-1})^s$ and a correction term, 
because the formulas then become much simpler than without this renormalization:

\renewcommand{\thetable}{\!\!}
\begin{table}[H]
\begin{tabular}{lll}
$n=3:$& $s=1:$& $1-L^2$\\
& $s=2:$&  $(1-L^2)^2$\\
& $s=3:$&  $(1-L^2)^3+144\,(1-5L^3+4L^6)$\\
& $s=4:$&  $(1-L^2)^4+576\,(1-94L^2-5L^3+245L^5+4L^6-151L^8)$\\
&&  $$\\
$n=4:$& $s=1:$& $1-L^3$\\
& $s=2:$& $(1-L^3)^2$\\
& $s=3:$& $(1-L^3)^3+\frac{36}{25}\,(4+72L-297L^5+221L^9)$\\
& $n=4$& $(1-L^3)^4+\frac{144}{125}\,(884+360L-20L^3-19584L^4-1485L^5$\\
&& \hspace{1.8in}$+ 44253L^8+1105L^9-25513L^{12})$\\
\\
$n=5:$& $s=1:$&  $1-L^4$\\
&  $s=2$&  $(1-L^4)^2$\\
& $s=3:$&  $(1-L^4)^3+\frac{32}{45}\,(7+134L^2-504L^7+363L^{12})$\\  
& $s=4:$&  $(1-L^4)^4+\frac{16}{135}\,(168+8576L+3216L^2-168L^4-127568L^6$\\
&& \hspace{1.4in}$-12096L^7+270144L^{11}+8712L^{12}-150984L^{16})$\\
\end{tabular}
\textit{
\caption{\normalfont List of values of
$\;s!\,\bigl(\frac{24n}{(n-2)(n+1)}\bigr)^s\Phi_s/L\;$  
for  $\,s=1,\,2,\,3,\,4\,$  and  $\,n=3,\,4,\,5$}}
\end{table}

\noindent This suggests that the series $\sum_s(\Phi_s/L)\,w^{-s}$ is given 
to a first approximation by a pure exponential $\,\exp\bigl(\frac{(n-2)(n+1)}{24n}(1-L^{n-1})/w\bigr)\,$ 
and hence that the formulas for the coefficients of the expansion~\e_ref{asym} 
may become simpler if we take the logarithm. Doing this, we find an expansion 
which begins
  $$ \log\F(w,x)\=\mu(x)w\+\log L(x)\+\frac{(n-2)(n+1)(1-L(x)^{n-1})}{24n}\,w\iv\+0\,w^{-2}\+\cdots $$
and in which, at least experimentally, the coefficient of $w^{-j}$ for $j\ge1$ 
is the sum of a term independent of~$x$
and a term of the form $L^{-j}$ times a polynomial (without constant term) in $L^n$.  
By applying the operator $\,w\iv D$ and adding 1, 
this can be stated more elegantly as follows.

{\it Conjecture:} 
If $\F$ is given by~\e_ref{F0def}, then
  \BE{LogExp}  1\+\frac xw\,\frac\p{\p x}\,\log\F(w,x)\;\overset?=\;L\,\sum_{k=0}^\infty\frac{P_k(n,L^n)}{(nLw)^k} \,,\EE
where $P_k(n,X)$ is a polynomial in $X$ of degree $k$ 
with coefficients in $\Q[n]$.

We have verified this conjecture up to order $\,\text O\big(w^{-6}\big)$,
with the values of the corresponding coefficients $P_k$ being given~by
\begin{alignat*}{1} 
P_0(n,X) &\= 1 \,,\\ 
P_1(n,X) &\= X-1\,,\\
P_2(n,X) &\= -\,\frac{(n+1)(n-1)(n-2)}{24}\,(X-1)X \,,  \\
P_3(n,X) &\= 0 \;,  \\
P_4(n,X) &\= \frac{(n+1)(n-1)(n-2)}{5760}\,(X-1)(A_3X^3+A_2X^2+A_1X)\,,  \\
P_5(n,X) &\= -\,\frac{(n+1)(n-1)(n-2)}{5760}\,(X-1)(B_4X^4+B_3X^3+B_2X^2+B_1X),
\end{alignat*}  
where
\begin{alignat*}{1} 
A_1 &\= (n-3)(7n^3-17n^2+22n-24)\,, \hspace{1.78in}\\
A_2 &\=  -\,(2n-3)(3n-1)(7n^2-9n+14)\,,\\
A_3 &\= 3\,(14n^4-33n^3+52n^2-23n-2)\,,\\
B_1 &\= -\,(n-3)(n-4)(7n^3-17n^2+22n-24)\,,\\ 
B_2 &\= 2\,(n-1)(n-2)(49n^3-115n^2+152n-124)\,,\\
B_3 &\= -4\,(n-1)(3n-1)(3n-4)(7n^2-9n+14)\,,\\
B_4 &\= 8\,(n-1)(3n-2)(7n^3-11n^2+17n-1)\,.
\end{alignat*}
The coefficients of the polynomials $P_k$ follow no apparent pattern apart from the divisibility by $(n+1)(n-1)(n-2)X(X-1)$:
the common factors of $A_1$ and $B_1$ and of $A_2$ and $B_3$ are striking, but nothing similar occurs for the next two polynomials.
On the other hand, there is a simple formula for the leading coefficient of $P_k(n,X)$ with respect to~$n$, namely (at least up to $k=7$) 
$$ P_k(n,X) \= \begin{cases} \qquad\a_j\,e_k(X)\,n^{4j-1}\+\O(n^{4j-2})& 
\text{if $k=2j>0$,}\\ 
(j-1)\a_j\,e_k(X)\,n^{4j}\+\O(n^{4j-1}) & \text{if $k=2j+1$}, \end{cases} $$
where $\a_j$ denotes the coefficient of $u^{2j}$ in $\frac{u/2}{\sinh u/2} $
($\a_0=1$, $\a_1=-\frac1{24}$, $\a_2=\frac7{5760}$, $\a_3=-\frac{31}{967680}$, \dots)
and where $e_1=X-1$,  $\,e_2=X^2-X$, $\,e_3=2X^3-3X^2+X$, \,\dots\ are the polynomials defined by
$$ e_k(X)\=\sum_{l=1}^k\,(-1)^{k-l}\,(l-1)!\,\S_k^{(l)}\,X^l \quad\in\quad\Z[X]$$
with $\S_k^{(l)}$ as before a Stirling number of the second kind.  
This is interesting because the argument $X=L^n$ of $P_k(n,X)$ in
equation~\e_ref{LogExp} is in fact $(1-n^nx)\iv$ and 
the functions $e_k((1-x)\iv)$ have the basic property
$$e_k\biggl(\frac1{1-x}\biggr)
\=\sum_{d=1}^\infty\,d^{k-1}x^d\quad\in\;x\,\Z[[x]]\qquad(k\ge1)\,. $$
There is also a possible intriguing connection with modular and elliptic functions since, for example,
the power series in two variables $\,\sum\a_je_{2j}\bigl(\frac1{1-x}\bigr)\,u^{2j-1}\,$ is closely related to the
expansion of the Weierstrass $\wp$-function and related Jacobi forms. 
This suggests possible
hidden modularity properties of the original function $\F(w,x)$.

As a final remark, we observe that~\e_ref{LogExp}, if it is true, defines the power series $\F(w,x)$ even for non-integral 
values of $n$ and shows that this function is analytic in~$n$ as well as in~$w$ and~$x$.  This seems surprising since~$\F$
is defined as a hypergeometric function of order~$n$ and we would usually not expect such series to have a reasonable interpolation
with respect to the order of the differential equation which they satisfy.

\vspace{.2in}

\noindent
{\it Max-Planck-Institut f$\ddot{u}$r Mathematik, Bonn\\
zagier@mpim-bonn.mpg.de}

\noindent
{\it Department of Mathematics, SUNY Stony Brook, NY 11794-3651\\
azinger@math.sunysb.edu}

\end{document}